\documentclass{article}
\usepackage{a4wide}
\usepackage{amsfonts}
\usepackage{amssymb}
\usepackage{amsthm}
\usepackage{enumerate}
\usepackage{amsmath}
\usepackage{times}
\usepackage{amscd}
\usepackage[all]{xy}
\usepackage{verbatim}
\usepackage{graphicx}
\usepackage{isolatin1}
\usepackage{hyperref}

\title{A polynomial with Galois group $\SL_2(\F_{16})$}
\author{Johan Bosman\footnote{Partially supported by the Dutch scientific organisation NWO.\ \ E-mail: jgbosman@math.leidenuniv.nl}}
\date{January 16, 2007}

\date\today

\theoremstyle{plain}
\newtheorem{theorem}{Theorem}
\newtheorem{proposition}{Proposition}
\newtheorem{conjecture}{Conjecture}

\theoremstyle{definition}

\newenvironment{proof2}[1]{\medskip \noindent {\bf #1}}
           {\unskip \nobreak \hfill \hbox{$\Box$} \par \medskip}

\newtheorem{exam}{Voorbeeld}

\renewcommand\mod{\mathop{\rm mod} \nolimits}
\newcommand\Frob{\mathop{\rm Frob} \nolimits}

\newcommand\F{\mathbb{F}}
\newcommand\Z{\mathbb{Z}}
\newcommand\Q{\mathbb{Q}}

\newcommand\C{\mathbb{C}}

\newcommand\Qbar{\overline{\mathbb{Q}}}
\renewcommand\P{\mathbb{P}}

\newcommand\el{\ell}

\newcommand\Tr{\mathop{\rm Tr} \nolimits}
\newcommand\Det{\mathop{\rm Det} \nolimits}
\newcommand\Aut{\mathop{\rm Aut} \nolimits}
\newcommand\GL{\mathop{\rm GL} \nolimits}
\newcommand\SL{\mathop{\rm SL} \nolimits}
\newcommand\Gal{\mathop{\rm Gal} \nolimits}

\renewcommand\ker{\mathop{\rm Ker} \nolimits}
\newcommand\im{\mathop{\rm Im} \nolimits}

\newcommand\tbox[1]{\mathop{\rm #1} \nolimits}

\newcommand\disc{\mathop{\rm Disc} \nolimits}
\renewcommand\O{\mathcal{O}}
\renewcommand\o[1]{\overline{#1}}

\newcommand\mat[4]{
  \left(
  {#1 \atop #3}
  \thinspace
  {#2 \atop #4}
  \right)
}

\makeatletter
\def\eqalign#1{\null\,\vcenter{\openup\jot\m@th
  \ialign{\strut\hfil$\displaystyle{##}$&$\displaystyle{{}##}$\hfil
      \crcr#1\crcr}}\,}

\newcounter{Lcount}
{%
  \medskip
  \setcounter{Lcount}{1}%
  \begin{list}{\bf\arabic{section}.\arabic{Lcount}.}{\usecounter{Lcount}%
    \settowidth\leftmargin{\bf\arabic{section}.99. }%
    \setlength\rightmargin{0pt}%
    \setlength\itemsep{\bigskipamount}}%
}%
{\end{list}}

{%
  \par
  \setlength{\leftmargin}{#1}%
  \setlength{\rightmargin}{#2}%
  \advance\linewidth -\leftmargin
  \advance\linewidth -\rightmargin
  \advance\@totalleftmargin\leftmargin
  \@setpar{{\@@par}}%
  \parshape 1 \@totalleftmargin \linewidth
  \ignorespaces
}%
{\par}

\makeatother

\begin{document}
\maketitle
\begin{abstract}
\noindent
In this paper we show that the polynomial
$
x^{17} - 5x^{16} + 12x^{15} - 28x^{14} + 72x^{13} - 132x^{12} + 116x^{11} - 74x^9 + 90x^8 - 28x^7 - 
12x^6 + 24x^5 - 12x^4 - 4x^3 - 3x - 1\in\Q[x]
$
has Galois group $\SL_2(\F_{16})$, filling in a gap in the tables of J\"urgen Kl\"uners and Gunther Malle (see \cite{KlMa}). The computation
of this polynomial uses modular forms and their Galois representations.
\end{abstract}
\section{Introduction}
It is a computational challenge to construct polynomials with a prescribed Galois group, see \cite{KlMa} for methods and examples. 
Here, by the Galois group of a polynomial
$f\in\Q[x]$ we mean the Galois group of a splitting field of $f$ over $\Q$
together with its natural action on the roots of $f$ in this splitting field.   
J\"urgen Kl\"uners informed me about an interesting group for which a polynomial had not been found yet, namely $\SL_2(\F_{16})$ with its natural
action on $\P^1(\F_{16})$. This action is faithful because of $\tbox{char}(\F_{16})=2$. It must be noted that the existence
of such a polynomial was already known to Mestre (unpublished). In this paper we will give an explicit example. 
\begin{proposition}\label{polynomial}
The polynomial
$$
\eqalign{
P(x) :=\,\, &x^{17} - 5x^{16} + 12x^{15} - 28x^{14} + 72x^{13} - 132x^{12} + 116x^{11}\cr& - 74x^9 + 90x^8 - 28x^7 - 
12x^6 + 24x^5 - 12x^4 - 4x^3 - 3x - 1\in\Q[x]
}
$$
has Galois group isomorphic to $\SL_2(\F_{16})$ with its natural action on $\P^1(\F_{16})$.
\end{proposition}\noindent
What is still unknown is
whether there exists a regular extension of $\Q(T)$ with Galois group isomorphic to $\SL_2(\F_{16})$; regular here means that it contains no
algebraic elements over $\Q$ apart from $\Q$ itself.
In section \ref{seccomput} we will say some words about the calculation of the polynomial and the connection with modular forms.
We'll indicate how one can
verify that it has the claimed Galois group in section
\ref{secverify} using computational Galois theory.
We will show in section \ref{secshowmod} that this polynomial gives a Galois representation associated to an explicitly given
modular form.

\section{Computation of the polynomial}\label{seccomput}
In this section we will briefly indicate how one can find a polynomial as in proposition \ref{polynomial}. We will make use of modular forms. 
For an overview as well as many further references on this subject the reader is refered to \cite{DiaIm}. 
\newline\smallskip

\noindent
Let $N$ be a positive integer and consider the space $S_2(\Gamma_0(N))$ of holomorphic cusp forms of
weight 2 for $\Gamma_0(N)$. A newform $f\in S_2(\Gamma_0(N))$ has a $q$-expansion $f=\sum_{n\geq 1} a_nq^n$
where the coefficients $a_n$ are in a number field. The smallest number field containing all the coefficients
is denoted by $K_f$. 
To a given prime number $\el$ and a place
$\lambda$ of $K_f$ above $\el$ one can attach a
semi-simple Galois representation $\o{\rho}_f:\Gal(\Qbar/\Q)\to\GL_2(\F_\lambda)$ unramified outside $N\el$
satisfying the following property: for each prime $p\nmid N\el$ and any Frobenius element
$\Frob_p\in\Gal(\o{\Q}/\Q)$ attached to $p$ we have 
\begin{equation}\label{trdet}
\Tr(\o{\rho}_f(\Frob_p))\equiv a_p\mod\lambda
\quad\mbox{and}\quad
\Det(\o{\rho}_f(\Frob_p))\equiv p\mod\lambda.
\end{equation}
The representation $\o{\rho}_f$ is unique up to isomorphism. 
The fixed field of $\ker(\o{\rho}_f)$ in $\Qbar$ is Galois over $\Q$ with Galois group
isomorphic to $\im(\o{\rho}_f)$. For $\el=2$ and any $\lambda$ above $\el$ equation (\ref{trdet}) together
with Chebotarev's density theorem
imply that $\im(\o{\rho}_f)$
is contained in $\SL_2(\F_\lambda)$. So to show that there is an extension of $\Q$ with Galois group
isomorphic to $\SL_2(\F_{16})$ it suffices to find an $N$ and a newform $f\in S_2(\Gamma_0(N))$ such that 
there is a prime $\lambda$ of degree 4 above 2 in $K_f$ and $\im(\o{\rho}_f)$ is the full group $\SL_2(\F_\lambda)$. 
Using modular symbols we can calculate the coefficients of
$f$, hence traces of matrices that occur in the image of $\o{\rho}_f$. 
For a survey paper on how this works, see \cite{St}.
A subgroup $\Gamma$ of $\SL_2(\F_{16})$ contains
elements of every trace if and only if $\Gamma$ equals $\SL_2(\F_{16})$; this can be shown in several ways, either by a direct calculation
or by invoking a more general classification result like \cite[Thm. III.6.25]{Su}. With this in mind, after a small computer search in which we check the occurring
values of $\Tr(\o{\rho}_f(\Frob_p))$ up to some moderate bound of~$p$,
one finds that a suitable modular form $f$ exists in $S_2(\Gamma_0(137))$. 
It turns out that we have $K_f\cong\Q(\alpha)$ with the minimal polynomial of $\alpha$ equal to $x^4+3x^3-4x-1$
and that $f$ is the form whose $q$-expansion starts with
$$
f = q + \alpha q^2 + (\alpha^3+\alpha^2-3\alpha-2)q^3 + (\alpha^2-2)q^4 + \cdots.
$$
\smallskip

\noindent
Now the next question comes in: knowing this modular form, how does one produce a polynomial? 
In general, one can use the Jacobian $J_0(N)$ to construct $\o{\rho}_f$. In this particular case we can do that in the following way.
We observe that $K_f$ is of degree 4 and that the prime 2 is inert in it. Furthermore we can verify that the subspace of $S_2(\Gamma_0(137))$ fixed by the Atkin-Lehner operator
$w_{137}$ is exactly the subspace generated by all the complex conjugates of $f$. 
These observations imply that $\o{\rho}_f$ is isomorphic to the action of $\Gal(\Qbar/\Q)$ on 
$\tbox{Jac}(X_0(137)/\langle w_{137}\rangle)[2]$, where we give this latter space an $\F_{16}$-vector space structure via the action of the Hecke operators. 
Note that $\im(\o{\rho}_f)=\SL_2(\F_{16})$
implies surjectivity of the natural map $\mathbb{T}\to\O_{K,f}/(2)\cong\F_{16}$, where $\mathbb{T}$ is the Hecke algebra attached to $S_2(\Gamma_0(N))$. 
The methods described in \cite{ECJ} allow us now to give complex approximations of 
the 2-torsion points of $\tbox{Jac}(X_0(137)/\langle w_{137}\rangle)$ to a high precision. 
This part of the calculation took by far the most
effort; the author will write more details about how this works in a future paper (or thesis).  
We use this to 
give a real approximation of a polynomial with Galois group isomorphic to $\SL_2(\F_{16})$. 
The paper \cite{ECJ} does, at least implicitly, give a theoretical upper bound for the height of the coefficients
of the polynomial hence an upper bound for the calculation precision to get an exact result. Though this upper
bound is small in the sense that it leads to a polynomial time algorithm, it is still far too high to be of use
in practice. However it turns out that we can use a much smaller precision to obtain our polynomial, the only
drawback being that this does not give us a proof of its correctness, so we have to verify this afterwards.
\newline\smallskip

\noindent
The polynomial $P'$ obtained in this way has coefficients of about 200 digits so we want to find a polynomial of smaller height defining the same number field
$K$. To do this, we first compute the ring of integers $\O_K$ of $K$. In \cite{BuLe} an algorithm to do this is described, provided that
one knows the squarefree factorisation of $\disc(f)$ and even if we don't know the squafrefree factorisation of the discriminant, the algorithm produces a
'good' order in $K$. 
Assuming that our polynomial $P'$ is correct we know that $K$ is unramified outside
$2\cdot 137$ so we can easily calculate the squarefree factorisation of $\disc(f)$ and hence apply the algorithm. 
Having done this we obtain an order in $K$ with a discriminant small enough to be able to factor and hence we know that this is indeed the maximal order $\O_K$.
Explicitly, the discriminant is equal to
\begin{equation}\label{discOK}
\disc(\O_K)=2^{30}\cdot 137^8.
\end{equation}
We embed $\O_K$ as a lattice into $\C^{[K:\Q]}$ in the natural way and use lattice basis reduction, see \cite{LLL}, to compute 
a short vector $\alpha\in\O_K-\Z$. The minimal polynomial of $\alpha$ has small coefficients. 
In our particular case $[K:\Q]$ is equal to 17, which is a prime number, hence this new polynomial must define the full field~$K$.
This method gives us also a way of expressing $\alpha$ as an element of $\Q(x)/(P'(x))$.

\section{Verification of the Galois group}\label{secverify}
Now that we have computed a polynomial $P(x)$,
we want to verify that its Galois group $\Gal(P)$ is really isomorphic to $\SL_2(\F_{16})$ and that we can identify the set $\Omega(P)$ of roots of $P$ with $\P^1(\F_{16})$ in such a way
that the action of $\Gal(P)$ on $\Omega(P)$ is identified with the action of $\SL_2(\F_{16})$ on $\P^1(\F_{16})$. 
\newline\smallskip

\noindent
For completeness let us remark that it is easy to verify that $P(x)$ is irreducible since it is irreducible modulo 5. 
The irreducibility of $P$ implies that $\Gal(P)$ is a transitive permutation group of degree 17. The transitive permutation groups of degree 17
have been classified, see for example \cite{Si}.
From \cite[Thm. III.6.25]{Su} it follows that up to conjugacy there is only one subgroup of index 17 in $\SL_2(\F_{16})$, namely the group of upper triangular
matrices. This implies that up to conjugacy there is exactly one transitive 
$G<S_{17}$ that is isomorphic to $\SL_2(\F_{16})$. 
Hence if $\Gal(P)\cong\SL_2(\F_{16})$ is an isomorphism of groups then there is an identification
of $\Omega(P)$ with $\P^1(\F_{16})$ such that the group actions become compatible. 
It follows from the classification in \cite{Si} that if the order of a transitive $G<S_{17}$ is divisible by 5, then $G$ must contain a 
transitive subgroup isomorphic to $\SL_2(\F_{16})$. To show that $5\mid \#\Gal(P)$ we use the fact that for a prime 
$p\nmid\disc(P)$ we the decomposition type of $P$ modulo $p$ is equal to the cycle type of any $\Frob_p\in\Gal(P)$ attached to $p$.
One can verify that modulo $7$ 
we get the following decomposition into irreducibles:
$$
\o{P}=
(x-3)(x-5)(x^{15} + 3x^{14} + 4x^{12} + 6x^{11} + 3x^{10} + x^9 + 5x^8 + x^6 + x^5 + 2x^4 + 4x^3 + 4x^2 + 3x + 6),
$$
showing that indeed $5\mid\#\Gal(P)$ hence $\Gal(P)$ contains $\SL_2(\F_{16})$ as a subgroup.
\newline\smallskip

\noindent 
To show that $\Gal(P)$ cannot be bigger than $\SL_2(\F_{16})$ it seems
inevitable to use heavy computer calculations. It would be interesting to see a method which does not use this.
\newline\smallskip

\noindent
Note that the action of $\SL_2(\F_{16})$ on $\P^1(\F_{16})$ is sharply 3-transitive. 
So first we show that $\Gal(P)$ is not 4-transitive to prove that it does not contain $A_{17}$. 
To do this we start with calculating the polynomial
\begin{equation}\label{Qprod}
Q(x) := \prod_{\{\alpha_1,\alpha_2,\alpha_3,\alpha_4\}\subset\Omega(P)} 
\left(X-\alpha_1-\alpha_2-\alpha_3-\alpha_4\right),
\end{equation}
where the product runs over all subsets of $\{1,\ldots,17\}$ consisting of exactly 4 elements.
This implies that $\deg(Q)=2380$. 
One can calculate $Q(x)$ by symbolic methods, see \cite{CaMc}.
Suppose that $\Gal(P)$ acting on
$\Omega(P)$ is 4-transitive. Then the action on $\Omega(Q)$ is transitive hence $Q(x)$ is irreducible.
So if we can show that $Q(x)$ is reducible, we have shown that $\Gal(P)$ is not 4-transitive. 
\newline\smallskip

\noindent
We have two ways to find a nontrivial factor of $Q(x)$: the first way is use a factorisation algorithm and the
second way is to produce a candidate factor ourselves. An algorithm that works very well for our
type of polynomial is Van Hoeij's algorithm, see \cite{Ho}. One finds that $Q(x)$ is the product of
3 polynomials of degrees 340, 1020 and 1020 respectively. A more direct way to produce a candidate
factorisation is as follows. The method from \cite{ECJ} gives a bijection between the set of approximated
complex roots
of $P'$ and the set $\P^1(\F_{16})$ such that the action of $\Gal(P')$ on $\Omega(P')$ 
corresponds to the action of $\SL_2(\F_{16})$ on $\P^1(\F_{16})$, assuming the outcome is correct.
From the previous section we know how to express the roots of $P$ as rational expressions in the
roots of $P$ hence this gives us a bijection between $\Omega(P)$ and $\P^1(\F_{16})$, conjecturally
compatible with the group actions of $\Gal(P)$ and $\SL_2(\F_{16})$ respectively. 
A calculation shows that the action of $\SL_2(\F_{16})$ on the set of 
unordered four-tuples of elements of
$\P^1(\F_{16})$ has 3 orbits, of size 340, 1020 and 1020 respectively. Using approximations to a
high precision of the roots, we use these orbits to produce sub-products of (\ref{Qprod}), round off
the coefficients to the nearest integer and verify afterwards that the obtained polynomials are
indeed factors of $Q(x)$. 
\newline\smallskip

\noindent
Let us remark that $\SL_2(\F_{16})\rtimes\Aut(\F_{16})$ with its natural action on $\P^1(\F_{16})$ is a transitive permutation group of degree
17, hence also its subgroup $G:=\SL_2(\F_{16})\rtimes\langle\Frob_2^2\rangle$. 
Furthermore, it is well-known that $\SL_2(\F_{16})\rtimes\Aut(\F_{16})$ is isomorpic to $\Aut(\SL_2(\F_{16}))$ (where $\SL_2(\F_{16})$ acts by
conjugation and $\Aut(\F_{16})$ acts on matrix entries) and actually inside $S_{17}$
this group is the normaliser of both $\SL_2(\F_{16})$ and itself.
According to the classifiation of
transitive permutation groups of degree 17 in \cite{Si} these two groups are the only ones that lie strictly between $\SL_2(\F_{16})$ and
$A_{17}$. Once we have fixed $\SL_2(\F_{16})$ inside $S_{17}$, these two groups are actually unique subgroups of $S_{17}$, not just up to
conjugacy.
\newline\smallskip

\noindent
Note that the index $[A_{17}:\Aut(\SL_2(\F_{16}))]$ is huge, namely $10897286400$.
If we can verify that $G$ is not a subgroup of $\Gal(P)$, then we are done. 
The fact that the index is small and the uniqueness of $G$ make an algorithm of Geissler and Kl\"uners,
see \cite{GeKl}, very
suitable to decide whether $G<\Gal(P)$. It turns out that this is not the case,
hence $\Gal(P)\cong \SL_2(\F_{16})$. 

\section{Does $P$ indeed define $\o\rho_f$?}\label{secshowmod}
So now that we have shown that $\Gal(P)\cong\SL_2(\F_{16})$ we can wonder whether we can prove that $P$ comes from the
modular form $f$ we used to construct it with. Once an isomorphism of $\Gal(P)$ with $\SL_2(\F_{16})$ is given,
the polynomial $P$ defines a representation $\o\rho_P:\Gal(\Qbar/\Q)\to\SL_2(\F_{16})$. 
Above we mentioned that 
that $\tbox{Out}(\SL_2(\F_{16}))$ is isomorphic to $\tbox{Aut}(\F_{16})$ acting on matrix
entries. Hence, up to an automorphism of $\F_{16}$, the map sending $\sigma\in\Gal(\Qbar/\Q)$ to the characteristic polynomial of
$\o\rho_P$ in $\F_{16}[x]$ is determined by $P$ and in fact the isomorphism class of $\o\rho_P$ is well-defined up to an automorphism of
$\F_{16}$.
More concretely, we have to show that the splitting field of $P$, which we will denote by $L$, 
is the fixed field
of $\ker(\o{\rho}_f)$.
In this section we will be using basic properties of local fields as can be found in
\cite{SeLF}.
\newline\smallskip

\noindent
A continuous representation $\o\rho:\Gal(\Qbar/\Q)\to\GL_2(\o\F_\el)$
has a \emph{level} $N(\o\rho)$ and a
\emph{weight} $k(\o\rho)$. Instead of repeating the full definitions here, which 
are lengthy (at least for the weight) and can be found in \cite{Se} 
(see also \cite{Ed}
for a discussion on the definition of the weight), we will just say that
they are defined in terms of the local representations $\o\rho_p:\Gal(\o{\Q}_p/\Q_p)\to\GL_2(\o\F_\el)$ obtained
from $\o\rho$. The level is defined in terms of the representations
$\o\rho_p$ with $p\ne\el$ and the weight is defined in terms of $\o\rho_\el$. 
Serre states the following conjecture in \cite{Se}.
\begin{conjecture}\label{serreconj}
Let $\el$ be a prime and let $\o\rho:\Gal(\Qbar/\Q)\to\GL_2(\o\F_\el)$ 
be a continuous odd irreducible Galois representation (a representation is called odd if the image of a complex conjugation has determinant
$-1$). Then there exists a modular form $f$ of level $N(\o\rho)$ and weight
$k(\o\rho)$ which is a normalised eigenform and a prime $\lambda\mid\el$ of $K_f$  
such that $\o\rho$ and $\o\rho_{f,\lambda}$ become isomorphic after a suitable embedding of $\F_\lambda$ into $\o\F_\el$.
\end{conjecture}\noindent
Recently, Khare and Wintenberger proved in \cite{KhWi} the following part of conjecture \ref{serreconj}.
\begin{theorem}\label{KhareWint}
Conjecture \ref{serreconj} holds in each of the following cases:
\begin{itemize}
\item $N(\o\rho)$ is odd and $\el>2$.
\item $\el=2$ and $k(\o\rho)=2$.
\end{itemize}
\end{theorem}\noindent
With theorem \ref{KhareWint} in mind it is sufficient to prove that a representation $\o\rho=\o\rho_P$ 
attached to $P$ has level $137$ and weight $2$, which
are the level and weight of the modular form $f$ we used to construct it with and that of all eigenforms in
$S_2(\Gamma_1(137))$, the form $f$ is one which gives rise to $\o\rho_P$. 
Therefore, in the remainder of this section we will verify the following proposition.
\begin{proposition}
Let $f$ be the cusp form from section \ref{seccomput}. Up to an automorphism of $\F_{16}$, the representations $\o\rho_P$ and
$\o\rho_{f,(2)}$ are isomorphic. In particular, the representation $\o\rho_P$ has Serre-level $137$ and Serre-weight $2$. 
\end{proposition}\noindent

\subsection{Verification of the level}
The level is the easiest of the two to verify. Here we have to do local computations in $p$-adic fields with $p\not=2$. According to the
definition of $N(\o\rho)$ in \cite{Se} it suffices to verify that $\o\rho$ is unramified outside 2 and 137, tamely ramified at $137$ 
and the local inertia subgroup $I$ at $137$ leaves exactly one
point of $\P^1(\F_{16})$ fixed. That $\o\rho_P$ is unramified outside 2 and 137 follows immediately from
(\ref{discOK}).
\newline\smallskip

\noindent
From (\ref{discOK}) and the fact that $137^8\Vert\disc(P)$ it follows that
the monogeneous order defined by $P$ is maximal at $137$. 
Modulo 137, the polynomial $P$ factors as
$$
\o{P}=(x+14)(x^2 + 6x + 101)^2(x^2 + 88x + 97)^2(x^2 + 106x + 112)^2(x^2 + 133x + 110)^2.
$$
Let $v$ be any prime above 137 in $L$. 
From the above factorisation it follows that the prime 137 decomposes in 
$K$ as a product of 5 primes; one of them has its inertial and ramification degree equal to 1 and the other four ones have their inertial and ramification degrees
equal to 2. 
Thus $\deg(v)$ is a power of 2, as $L$ is obtained by succesively adjoining roots of $P$ and in each step the relative
inertial and ramification degrees of the prime below $v$ are both at most 2. In particular, $\Gal(L_v/\Q_{137})$ is a subgroup of
$\SL_2(\F_{16})$ whose order is a power of 2. Now, $\{\mat{1}{*}{0}{1}\}$ is a Sylow 2-subgroup of $\SL_2(\F_{16})$, so $\Gal(L_v/\Q_{137})$ is, up to conjugacy,
a subgroup of $\{\mat{1}{*}{0}{1}\}$. Hence $I$ is also conjugate to a subgroup of $\{\mat{1}{*}{0}{1}\}$ and it is actually nontrivial
because 137 ramifies in $L$ (so $I$ of order 2 since the tame inertia group of any finite Galois extension of local fields
is cyclic).
\newline\smallskip

\noindent
It is immediate that $\o\rho$ is tamely ramified at 137 as no power of 2 is divisible by 137. Also, it is clear that $I$ has exactly one fixed
point in $\P^1(\F_{16})$ since $[{*\choose 0}]$ is the only fixed point of any nontrivial element of $\{\mat{1}{*}{0}{1}\}$. This establishes
the verification of $N(\o\rho)=137$.

\subsection{Verification of the weight}
Because the weight is defined in terms of the induced local representation $\Gal(\Qbar_2/\Q_2)$, we will try to compute some relevant properties of the
splitting field $L_v$ of $P$ over $\Q_2$, where $v$ is any place of $L$ above $2$. 
In $p$-adic fields one can only do calculations with a certain precision, but this does not give any problems since practically all properties one needs to
know can be verified rigorously using a bounded precision calculation and the error bounds in the calculations can be kept track of exactly.  
\newline\smallskip

\noindent
The polynomial $P$ does not define an order which is maximal at the prime 2. Instead we use the polynomial
$$
\eqalign{R&=
x^{17} - 11x^{16} + 64x^{15} - 322x^{14} + 916x^{13} + 276x^{12} - 5380x^{11} + 2748x^{10} + 6904x^9
 - 23320x^8\cr& + 131500x^7 - 140744x^6 -
16288x^5 - 39752x^4 - 48840x^3 + 102352x^2 + 234466x - 1518,
}
$$
which is the minimal polynomial of 
$$
\eqalign{
&\big(36863 + 22144\alpha + 123236\alpha^2 + 154875\alpha^3 - 416913\alpha^4 + 436074\alpha^5 + 229905\alpha^6 -
1698406\alpha^7\cr& + 1857625\alpha^8 - 467748\alpha^9 - 2289954\alpha^{10} + 2838473\alpha^{11} - 1565993\alpha^{12} +
605054\alpha^{13} - 263133\alpha^{14}\cr& + 112104\alpha^{15} - 22586\alpha^{16}\big)/8844,
}
$$
where $\alpha$ is a root of $P$. 
We can factor $R$ over $\Q_2$ and see that it has one root in $\Q_2$ which happens to be odd, and an Eisenstein factor of degree 16, which we will
call $E$. This type of decomposition can be read off from the Newton polygon of $R$ and it also shows that the order defined by $R$ is indeed maximal at 2.
From the oddness of the root and (\ref{discOK}) it follows that.
\begin{equation}\label{discE}
v_2(\disc(E))=30.
\end{equation}
For the action of $\Gal(\o{\Q}_2/\Q_2)$ on $\P^1(\F_{16})$ the factorisation means that there is one fixed point and one orbit
of degree 16. 
If we adjoin a root $\beta$ of $E$ to $\Q_2$ and factor $E$ over $\Q_2(\beta)$ then we see that it has an
irreducible factor of degree 15; 
in \cite{CaGo} one can find 
methods for factorisation and irreducibility testing that can be used to verify this. This means that $[L_v:\Q_2]$ is at least 240.
\newline\smallskip

\noindent
A subgroup of $\SL_2(\F_{16})$ that fixes a point has to be conjugate to a subgroup of the group
$$
H := \left\{\mat{*}{*}{0}{*}\right\}\subset\SL_2(\F_{16}),
$$
which is the stabiliser subgroup of $[{*\choose 0}]$. From $\#H=240$ it follows that $\Gal(L_v/\Q_2)$ 
is isomorphic to $H$ and from now on we will identify these two groups with each other. 
We can filter $H$ by normal subgroups:
$$
H\supset I \supset I_2\supset\{e\},
$$
where $I$ is the inertia subgroup and $I_2$ is the wild ramification subgroup, 
which is the unique Sylow 2-subgroup of $I$. 
We wish to determine the groups $I$ and $I_2$. Let $k(v)$ be the residue class field
of $L_v$. The group $H/I$ is isomorphic to $\Gal(k(v)/\F_2)$
and $I/I_2$ is isomorphic to a subgroup of $k(v)^*$. In particular
$[I:I_2]\mid (2^{[H:I]}-1)$ follows. 
The group $H$ has the nice property 
$$
[H,H]=\left\{\mat{1}{*}{0}{1}\right\}\cong \F_{16},
$$
which is its unique Sylow 2-subgroup. As $H/I$ is abelian, we see that $[H,H]\subset I$. 
We conclude that $I_2=[H,H]$, since above we remarked that
$I_2$ is the unique Sylow $2$-subgroup of $I$. The restriction $[I:I_2]\mid (2^{[H:I]}-1)$ leaves only one
possibility for $I$, namely $I=I_2$. 
\newline\smallskip

\noindent
Let $L_v'$ be the subextension of $L_v/\Q_2$ fixed by $I$. Then $L_v'$ is the maximal unramified 
subextension as well as the maximal tamely ramified subextension. It is in fact isomorphic to $\Q_{2^{15}}$,
the unique unramified
extension of $\Q_2$ of degree 15 and the Eisenstein polynomial $E$ from above, being irreducible over any
unramified extension of $\Q_2$, is a defining polynomial for the
extension $L_v/\Q_{2^{15}}$. 
According to \cite[Thm.~3]{MoTa} we can relate the discriminant of $L_v$ to $k(\o\rho)$ as follows:
$$
v_2(\disc(L_v)) = 
\left\{\begin{array}{ll}
240\cdot\frac{15}{8}=450 &\mbox{if $k(\o\rho)=2$}\cr
240\cdot\frac{19}{8}=570 &\mbox{if $k(\o\rho)\not=2$}
\end{array}\right.
$$
It follows from (\ref{discE}) that $v_2(\disc(L_v/\Q_2)) = 30\cdot 15 = 450$, so indeed $k(\o\rho)=2$.

\subsection{Verification of the form $f$}
Now we know $N(\o\rho_P)=137$ and $k(\o\rho_P)=2$ theorem \ref{KhareWint}
shows that there is an eigenform $g\in S_2(\Gamma_1(137))$ giving rise to $\rho_P$.
Using \cite[Cor. 2.7]{Bu} we see that if such a
$g$ exists, then there actually exists such a $g$ of trivial Nebentypus, i.e. $g\in S_2(\Gamma_0(137))$ (as $\SL_2(\F_{16})$ is
non-solvable $\rho_P$ cannot be an induced Hecke character from $\Q(i)$).
\newline\smallskip

\noindent
A modular symbols calculation shows that there exist two Galois orbits of newforms in $S_2(\Gamma_0(137))$: the form $f$ we used for our calculations and
another form, $g$ say. The prime 2 decomposes in $K_g$ as a product $\lambda^3\mu$, where $\lambda$ has inertial degree 1 and $\mu$ has
inertial degree $4$. So it could be that $g\mod\mu$ gives rise to $\rho_P$. We will show now that $f\mod\, (2)$ and $g\mod\mu$ 
actually give the same
representation. The completions of $\O_{K_f}$ and $\O_{K_g}$ at the primes $(2)$
and $\mu$ respectively are both isomorphic to $\Z_{16}$, the unramified extension of $\Z_2$ of degree 4. After a choice of embeddings of
$\O_{K_f}$ and $\O_{K_g}$ into $\Z_{16}$ we obtain two modular forms $f'$ and $g'$ with coefficients in $\Z_{16}$ and we wish to show that a suitable
choice of embeddings exists such that they are
congruent modulo 2. According to \cite[Thm. 1]{Stu}, it suffices to check there is a suitable choice of embeddings that gives
$a_n(f')\equiv a_n(g')\mod 2$ for all 
$n\leq [\SL_2(\Z):\Gamma_0(137)]/6=23$ (in \cite{Stu} this theorem is formulated for modular forms with coefficients in the ring of integers
of a number field, but the proof also works for $p$-adic rings). Using a modular symbols calculation, this can be easily verified.

\section{Acknowledgements}
I would like to thank J\"urgen Kl\"uners for proposing this computational challenge and explaining some computational Galois theory to me. Furthermore I 
want to thank Bas Edixhoven for teaching me about modular forms and the calculation of their coefficients. 
All the calculations were done with MAGMA (see \cite{Magma}), many of them on the MEDICIS cluster (\url{http://medicis.polytechnique.fr}). For being able to make use of the cluster I want to thank Marc Giusti
and Pierre Lafon.

\end{document}